\documentclass[a4paper, 12pt] {article}
\usepackage[cp1251]{inputenc}
\usepackage[russian, english]{babel}
\usepackage{amsfonts}
\usepackage{mathtext}
\usepackage{cite}

\usepackage{graphicx}

\usepackage{verbatim}
\usepackage{amsmath}
\usepackage{amsthm}
\usepackage{amssymb}
\usepackage{delarray}
\usepackage{mathrsfs}
\usepackage{xcolor}

\begin{document}

\thispagestyle{empty}

\begin{center}
{\Large\bf On recovering Sturm--Liouville-type operators with global delay on graphs from two spectra }

\end{center}

\begin{center}
{\large\bf Sergey Buterin\footnote{Department of Mathematics, Saratov State University, Russia {\it email: buterinsa@sgu.ru}} }
\end{center}

{\bf Abstract.} We suggest a new formulation of the inverse spectral problem for second-order functional-differential operators on
star-shaped graphs with global delay. The latter means that the delay, being measured in the direction to a specific boundary vertex, called
the root, propagates through the internal vertex to other edges. Now, we intend to recover the potentials given the spectra of two boundary
value problems on the graph with a common set of boundary conditions at all boundary vertices except the root. We prove the uniqueness
theorem and obtain a constructive procedure for solving this inverse problem assuming that the common boundary conditions are of the Robin
type and they are pairwise linearly independent. Although we focus on graphs with equal edges when the delay parameter also coincides with
their length, the proposed formulation is expected to be relevant for more general situations including non-star trees with non-equal edges
and with a wide range for the global delay parameter.


\\

{\large\bf 1. Introduction}
\\

In this paper, we test another formulation of the inverse spectral problem for Sturm--Liouville-type operators on star-shaped graphs with
global delay, which were introduced in \cite{But23}. This formulation considerably reduces the input data compared with those used for
recovering classical Sturm--Liouville operators on trees. Moreover, the new formulation allows to gather all necessary spectral information
by altering boundary conditions at only one boundary vertex, which becomes a new quality in the inverse spectral theory on graphs.

Purely differential (local) operators on geometrical graphs, often called quantum graphs, model various processes in networks and branching
structures appearing in organic chemistry, mesoscopic physics, quantum mechanics, nanotechnology, hydrodynamics, waveguide theory, and other
applications \cite{Mont, Nic, vB, LangLeug, Kuch, BCFK, BerkKuch, Pok, Kuz-17, Boris-22}. The theory of quantum graphs includes numerous
studies related to inverse problems of recovering operators from their spectral characteristics (see \cite{Ger, Bel, BrauWeik, Yur-05,
Pivo07, Ign15, Yurko-16, Bond-20, KravchAvd} and references therein). In particular, it was established in \cite{Yur-05} that for the unique
determination of the Sturm--Liouville potentials on all edges of an arbitrary compact tree with $m$ boundary vertices, it is sufficient to
specify the spectra of precisely $m$ specially chosen boundary value problems. The way how those spectra are gathered generalizes the
classical inverse Sturm--Liouville problem due to Borg \cite{B, novabook, BK19}. Moreover, such input data remain minimal among other
formulations of inverse spectral problems for purely differential operators on trees.

Concerning {\it functional}-differential as well as other classes of nonlocal operators on graphs, they were studied mostly in the {\it
locally} nonlocal case when the corresponding nonlocal equation on each edge can be treated independently of the other edges \cite{Nizh-12,
Bon18-1, HuBondShYan19, Hu20, WangYang-21, Bon22, WangYang-22}. This limitation has left uncertain how the nonlocalities could coexist with
the internal vertices of the graph.

In \cite{But23}, however, we introduced Sturm--Liouville-type operators on geometrical graphs with a global delay that naturally extends
through the internal vertices. On a star-shaped graph consisting of $m$ equal edges, this idea can be illustrated by the boundary value
problem
\begin{equation}\label{1.8}
-y_j''(x)+q_j(x)y_j(x-a)=\lambda y_j(x), \quad 0<x<1, \quad j=\overline{1,m},
\end{equation}
\begin{equation}\label{1.9}
y_j(x-a)=y_1(x-a+1), \quad \max\{0,a-1\}<x<\min\{a,1\}, \quad j=\overline{2,m},
\end{equation}
\begin{equation}\label{1.10}
y_1(1)=y_j(0), \quad j=\overline{2,m}, \qquad y'_1(1)=\sum_{j=2}^m y_j'(0),
\end{equation}
\begin{equation}\label{1.11}
y_1^{(\nu_1)}(0)+\nu_1 H_1y_1(0)=0,\quad V_j(y_j):=y_j^{(\nu_j)}(1)+\nu_j H_jy_j(1)=0, \quad j=\overline{2,m},
\end{equation}
where $a\in (0,2),$ while $\nu_j\in\{0,1\}$ and $H_j\in{\mathbb C}$ for $j=\overline{1,m}.$ It is also assumed that all potentials $q_j(x)$
belong to $L_2(0,1)$ and obey the conditions
\begin{equation}\label{1.11-2}
q_1(x)=0,\;\; x\in(0,\min\{a,1\}), \qquad q_j(x)=0, \;\; x\in(0,\max\{0,a-1\}), \;\; j=\overline{2,m},
\end{equation}
which make the equations in (\ref{1.8}) well defined. Supposing the $j$-th equation in (\ref{1.8}) to be defined on the edge $e_j$ of the
graph $\Gamma_m$ illustrated on Fig.~1., one can observe that relations (\ref{1.10}) specify the standard matching conditions at the internal
vertex $v_1,$ while (\ref{1.11}) become Dirichlet and/or Robin boundary conditions at the boundary vertices $v_0$ and $v_j$ for
$j=\overline{2,m}.$

\medskip
\begin{center}
\unitlength=0.7mm
\begin{picture}(80,90)
 \put(40,46){\line(2,3){27.7}}
 \put(40,46){\line(-1,0){50}}
 \put(40,46){\line(5,2){46.65}}
 \multiput(40,46)(2.5,-1){19}{\circle*{0.7}}
 \put(40,46){\line(2,-3){27.7}}

 \put(-10,46){\circle*{1}}
 \put(40,46){\circle*{1}}
 \put(68,88){\circle*{1}}
 \put(87,64.8){\circle*{1}}
 \put(68,4){\circle*{1}}

 \put (11,41){\small $e_1$}
 \put (55,67){\small $e_2$}
 \put (64,52.5){\small $e_3$}
 \put (47,21){\small $e_m$}

 \put (-18,44){\small $v_0$}

 \put (35,41){\small $v_1$}

 \put (68,90){\small $v_2$}
 \put (88,66){\small $v_3$}
 \put (68,0){\small $v_m$}

 \put (-10,48){\small $0$}
 \put (34,48){\small $1$}
 \put (38,51){\small $0$}
 \put (44,49.8){\small $0$}
 \put (42.34,42.0){\small $0$}
 \put (61,86){\small $1$}
 \put (82,66){\small $1$}
 \put (65,8){\small $1$}

 \put (10,-10){\small Fig. 1. Graph $\Gamma_m$}

\end{picture}
\end{center}

\bigskip
\bigskip
Relations (\ref{1.8})--(\ref{1.11-2}) generate a nonlocal operator with the global delay $a.$ The latter is measured in the direction to the
boundary vertex $v_0,$ which we refer to as {\it root}. Conditions (\ref{1.9}) specify an initial function for all equations in (\ref{1.8})
except the first one. In other words, those equations always involve the unknown function $y_1(x)$ on edge $e_1.$ This means that the delay
“extends” or “propagates” through the internal vertex $v_1.$ In \cite{But23}, this new quality was generalized in a natural way also to
operators on an arbitrary tree with variable edge lengths.

Moreover, in \cite{But23}, an inverse problem was studied for the particular case when $m=3$ and $a=1$ that consisted in recovering the
remaining functions $q_2(x)$ and $q_3(x)$ from the spectra corresponding to the two sets of parameters: $\nu_1=\nu_2=\nu_3-1=0$ and
$\nu_1=\nu_2-1=\nu_3=0,$ while $H_2=H_3=0.$ That statement resembles the inverse problem posed in \cite{Yur-05} for the classical
Sturm--Liouville operator $(a=0)$ on a tree, when the input spectral data are gathered by altering the boundary conditions at all boundary
vertices except the root.

In the present paper, we suggest another statement of the inverse problem that is assumed to involve always only two spectra obtained by
altering the boundary conditions namely at the root, while at the rest boundary vertices, the corresponding conditions remain fixed. This new
statement will be tested below for arbitrary $m$ and $a\in[1,2).$

For $\nu=0,1,$ denote by $\{\lambda_{n,\nu}\}_{n\in{\mathbb N}}$ the spectrum of the eigenvalue problem ${\cal B}_\nu:={\cal B}_\nu(q,H)$
consisting of relations (\ref{1.8})--(\ref{1.11-2}) under the settings:
$$
a\in[1,2), \quad \nu_1=\nu, \quad H_1=0, \quad \nu_j=1, \;\; j=\overline{2,m},
$$
and depending on $q:=[q_2,\ldots,q_m]$ and $H:=[H_2,\ldots,H_m].$ Indeed, according to (\ref{1.11-2}), admitted values of $a$ automatically
imply $q_1=0,$ while the remaining $q_j(x)$ vanish a.e. on $(0,a-1).$

Consider the following inverse problem.

\medskip
{\bf Inverse Problem 1.} Given $\{\lambda_{n,0}\}_{n\in{\mathbb N}}$ and $\{\lambda_{n,1}\}_{n\in{\mathbb N}}$ as well as the value of $a,$
find $q_j(x)$ for $j=\overline{2,m}$ provided that $H_j,\,j=\overline{2,m},$ are known a priori and distinct.

\medskip
We note that the distinction as well as the knowledge of the coefficients $H_j$ are, obviously, vital whenever one aims to recover each
$q_j(x)$ with reference to the corresponding edge $e_j.$

The following uniqueness theorem holds.

\medskip
{\bf Theorem 1. }{\it Let the coefficients $H_j,\,j=\overline{2,m},$ be known and distinct. Then the specification of
$\{\lambda_{n,0}\}_{n\in{\mathbb N}}$ and $\{\lambda_{n,1}\}_{n\in{\mathbb N}}$ as well as $a$ uniquely determines $q_j(x)$ for all
$j=\overline{2,m}.$}

\medskip
The case of smaller values of $a$ is more difficult and requires a separate investigation. However, an analogous uniqueness result is
expected to hold also for a certain range of $a<1$ and even for non-star trees with appropriately chosen the common boundary conditions.

Moreover, as for the inverse problem studied in \cite{But23}, one can expect that an analogous assertion will hold also after reducing the
input data to appropriate subspectra of the problems ${\cal B}_0$ and ${\cal B}_1,$ but any such refinement is beyond the goals of the
present paper.

Note that the problem ${\cal B}_\nu$ in the case $m=2$ becomes a usual Sturm--Liouville-type problem with constant delay on an interval:
\begin{equation}\label{1.11-3}
-y''(x)+q(x)y(x-a)=\lambda y(x), \quad 0<x<2, \qquad y^{(\nu)}(0)=y'(2)+H_2y(2)=0,
\end{equation}
where
$$
y(x)=\left\{\begin{array}{cc}y_1(x), & 0\le x\le 1,\\[3mm]
y_2(x-1), & 1< x\le 2,
\end{array}\right.
 \quad
q(x)=\left\{\begin{array}{cc}q_1(x), & 0< x<1,\\[3mm]
q_2(x-1), & 1< x< 2,
\end{array}\right.
$$
while $q_1(x)$ will actually vanish if our assumption about $a$ remains.

Although the inverse problems for functional-differential operators with constant delay on an interval have been studied quite extensively
(see, e.g., \cite{Pik91, FrYur12, Yang, Ign18, BondYur18-1, ButYur19, VPV19, DV19, SatShieh19, WShM19, Dur20, DB21, DB21-2, DB22, BMSh21,
ButDjur22, WKSh23, ButVas23, DjurVojv23}), Inverse Problem~1 in the case $m=2$ does not have any direct analog in the literature. Indeed, all
studies of recovering $q(x)$ in the equation (\ref{1.11-3}) address the case of a common boundary condition at the left-hand rather than the
right-hand end-point of the interval as in (\ref{1.11-3}). In other words, the spectra usually used for recovering $q(x)$ actually
corresponded to the boundary conditions
$$
y^{(\nu_1)}(0)+\nu_1 H_1y(0)=y^{(\nu)}(2)=0, \quad \nu=0,1,
$$
where $\nu_1\in\{0,1\}$ was fixed. While $q(x)$ is uniquely determined by those spectra whenever $a\in[4/5,2)$ (see \cite{BondYur18-1,
VPV19}), the recent series of papers \cite{DB21, DB21-2, DB22} establishes its possible non-uniqueness for any $a\in(0,4/5)$ (see also a
brief survey in \cite{BMSh21}). The present study, in particular, raises analogous questions for $a<1$ also in the case of the boundary
conditions in (\ref{1.11-3}).

The paper is organized as follows. In the next section, we construct the characteristic determinants of the problems ${\cal B}_\nu.$ The
proof of Theorem~1 along with a constructive procedure for solving Inverse Problem~1 is given in Section~3.
\\

{\large\bf 2. Characteristic determinants}
\\

Denote
\begin{equation}\label{02.0}
S_0(x,\lambda):=\frac{\sin\rho x}\rho, \quad S_1(x,\lambda):=\cos\rho x.
\end{equation}
For $\nu=0,1$ and $j=\overline{2,m},$ we also introduce the designations
\begin{equation}\label{02.1}
v_{\nu,j}(\lambda):=V_j(S_\nu(\,\cdot\,,\lambda)), \quad Q_{\nu,j}(x,\lambda)=\int_{a-1}^x \frac{\sin\rho(x-t)}\rho
q_j(t)S_\nu(t-a+1,\lambda)\,dt, \quad \rho^2=\lambda,
\end{equation}
where $V_j(\,\cdot\,)$ are the linear forms defined in (\ref{1.11}). Consider the entire functions
\begin{equation}\label{02.4}
\Delta_\nu^0(\lambda)=S_\nu'(1,\lambda)\prod_{j=2}^m v_{0,j}(\lambda) +S_\nu(1,\lambda)\sum_{j=2}^m v_{1,j}(\lambda)\prod_{{l\ne
j}\atop{l=2}}^m v_{0,l}(\lambda), \quad \nu=0,1,
\end{equation}
and
\begin{equation}\label{02.5}
\Delta_\nu(\lambda)=\Delta_\nu^0(\lambda) +\sum_{j=2}^m V_j(Q_{\nu,j}(\,\cdot\,,\lambda))\prod_{{l\ne j}\atop{l=2}}^m v_{0,l}(\lambda), \quad
\nu=0,1.
\end{equation}
Here and below, the prime always means the differentiation with respect to the {\it first} argument.

The function $\Delta_\nu(\lambda)$ is called {\it characteristic determinant} or {\it characteristic function} of the problem ${\cal B}_\nu.$
The following lemma holds.

\medskip
{\bf Lemma 1. }{\it For $\nu=0,1,$ eigenvalues of the problem ${\cal B}_\nu$ coincide with zeros of $\Delta_\nu(\lambda).$}

\medskip
{\it Proof.} In accordance with \cite{But23}, the functions
\begin{equation}\label{02.6}
y_1(x)= \sum_{l=0}^1C_{l,1}S_l(x,\lambda), \quad y_j(x)=\sum_{l=0}^1C_{l,j}S_l(x,\lambda) +\sum_{l=0}^1C_{l,1}Q_{l,j}(x,\lambda), \;\;
j=\overline{2,m},
\end{equation}
where $C_{l,j}$ are indefinite constants, form the so-called global general solution of the system of equations (\ref{1.8}) subject to the
initial-function conditions (\ref{1.9}). Substituting (\ref{02.6}) into the matching conditions (\ref{1.10}) as well as the boundary
conditions of ${\cal B}_\nu,$ we arrive at a system of linear algebraic equations with respect to the vector
$[C_{\nu,1},C_{0,2},C_{1,2},\ldots,C_{0,m},C_{1,m}]^T,$ having the determinant
$$
D_\nu(\lambda):=\left|\begin{array}{cccccccc}
S_\nu(1,\lambda)  &  0 & -1 &  0 &  0 & \ldots &  0 & 0 \\[3mm]
S_\nu(1,\lambda)  &  0 &  0 &  0 & -1 & \ldots &  0 & 0 \\[3mm]
\vdots & \vdots  & \vdots & \vdots &  \vdots &     & \vdots & \vdots \\[3mm]
S_\nu(1,\lambda)  &  0 &  0 &  0 &  0 & \ldots &  0 & -1 \\[3mm]
S_\nu'(1,\lambda) & -1 &  0 & -1 &  0 & \ldots & -1 &  0 \\[3mm]
V_2(Q_{\nu,2}(\,\cdot\,,\lambda)) & v_{0,2}(\lambda) &  v_{1,2}(\lambda) & 0 &  0 & \ldots & 0 & 0 \\[3mm]
V_3(Q_{\nu,3}(\,\cdot\,,\lambda)) & 0 &  0 & v_{0,3}(\lambda) &  v_{1,3}(\lambda) & \ldots & 0 & 0 \\[3mm]
\vdots & \vdots  & \vdots & \vdots &  \vdots &     & \vdots & \vdots \\[3mm]
V_m(Q_{\nu,m}(\,\cdot\,,\lambda)) & 0 &  0 & 0 & 0 & \ldots & v_{0,m}(\lambda) &  v_{1,m}(\lambda)
\end{array}\right|
$$
of order $2m-1$ (the $m$-th entry of the first column is $S_\nu'(1,\lambda)).$

Denote by $D_\nu^0(\lambda)$ the determinant that is obtained from $D_\nu(\lambda)$ after taking all $q_j(x)$ to equal zero, i.e. by zeroing
the last $m-1$ entries of the first column. Then we have
\begin{equation}\label{02.6-1}
D_\nu(\lambda)=D_\nu^0(\lambda) +\sum_{j=m+1}^{2m-1} V_{j-m+1}(Q_{\nu,{j-m+1}}(\,\cdot\,,\lambda))A_j(\lambda), \quad \nu=0,1,
\end{equation}
where $A_j(\lambda)$ is the cofactor of the $j$-th entry of the first column. In particular, we have
\begin{equation}\label{02.7}
A_m(\lambda)=(-1)^{m+1} \left|\begin{array}{ccccccc}
0 & -1 &  0 &  0 & \ldots &  0 & 0 \\[3mm]
0 &  0 &  0 & -1 & \ldots &  0 & 0 \\[3mm]
\vdots  & \vdots & \vdots &  \vdots &     & \vdots & \vdots \\[3mm]
0 &  0 &  0 &  0 & \ldots &  0 & -1 \\[3mm]
v_{0,2}(\lambda) &  0 & 0 &  0 & \ldots & 0 & 0 \\[3mm]
0 &  0 & v_{0,3}(\lambda) &  0 & \ldots & 0 & 0 \\[3mm]
\vdots  & \vdots & \vdots &  \vdots &     & \vdots & \vdots \\[3mm]
0 &  0 & 0 & 0 & \ldots & v_{0,m}(\lambda) &  0
\end{array}\right| =\sigma_m \prod_{l=2}^m v_{0,l}(\lambda),
\end{equation}
where $\sigma_m=(-1)^\frac{(m-1)(m+2)}2.$ For $j=\overline{m+1,2m-1},$ we have
$$
A_j(\lambda)=(-1)^{j+1} \left|\begin{array}{ccccccc}
0 & -1 &  0 &  0 & \ldots &  0 & 0 \\[3mm]
0 &  0 &  0 & -1 & \ldots &  0 & 0 \\[3mm]
\vdots  & \vdots & \vdots &  \vdots &     & \vdots & \vdots \\[3mm]
0 &  0 &  0 &  0 & \ldots &  0 & -1 \\[3mm]
 -1 &  0 & -1 &  0 & \ldots & -1 &  0 \\[3mm]
v_{0,2}(\lambda) &  0 & 0 &  0 & \ldots & 0 & 0 \\[3mm]
0 &  0 & v_{0,3}(\lambda) &  0 & \ldots & 0 & 0 \\[3mm]
\vdots  & \vdots & \vdots &  \vdots &     & \vdots & \vdots \\[3mm]
0 &  0 & 0 & 0 & \ldots & v_{0,m}(\lambda) &  0
\end{array}\right|,
$$
where the row possessing $v_{0,j-m+1}(\lambda)$ is absent. Thus, we obtain the representation
\begin{equation}\label{02.8}
A_j(\lambda)=\sigma_m \prod_{{l\ne j-m+1}\atop{l=2}}^m v_{0,l}(\lambda), \quad j=\overline{m+1,2m-1},
\end{equation}
which also fits (\ref{02.7}) for $j=m.$ Substituting (\ref{02.8}) into (\ref{02.6-1}) we get
\begin{equation}\label{02.9}
D_\nu(\lambda)=D_\nu^0(\lambda) +\sigma_m \sum_{j=m+1}^{2m-1} V_{j-m+1}(Q_{\nu,{j-m+1}}(\,\cdot\,,\lambda))\prod_{{l\ne j-m+1}\atop{l=2}}^m
v_{0,l}(\lambda), \quad \nu=0,1.
\end{equation}

For calculating $D_\nu^0(\lambda),$ it is convenient to write it in the form
$$
D_\nu^0(\lambda)=\left|\begin{array}{cccccccc}
S_\nu(1,\lambda)  &  0 & -1 &  0 &  0 & \ldots &  0 & 0 \\[3mm]
S_\nu(1,\lambda)  &  0 &  0 &  0 & -1 & \ldots &  0 & 0 \\[3mm]
\vdots & \vdots  & \vdots & \vdots &  \vdots &     & \vdots & \vdots \\[3mm]
S_\nu(1,\lambda)  &  0 &  0 &  0 &  0 & \ldots &  0 & -1 \\[3mm]
S_\nu'(1,\lambda) & -1 &  0 & -1 &  0 & \ldots & -1 &  0 \\[3mm]
v_{1,2}(\lambda) S_\nu(1,\lambda) & v_{0,2}(\lambda) &  0 & 0 &  0 & \ldots & 0 & 0 \\[3mm]
v_{1,3}(\lambda) S_\nu(1,\lambda) & 0 &  0 & v_{0,3}(\lambda) &  0 & \ldots & 0 & 0 \\[3mm]
\vdots & \vdots  & \vdots & \vdots &  \vdots &     & \vdots & \vdots \\[3mm]
v_{1,m}(\lambda)S_\nu(1,\lambda) & 0 &  0 & 0 & 0 & \ldots & v_{0,m}(\lambda) &  0
\end{array}\right|,
$$
which immediately implies
$$
D_\nu^0(\lambda) =S_\nu'(1,\lambda)A_m(\lambda) +S_\nu(1,\lambda)\sum_{j=m+1}^{2m-1} v_{1,j-m+1}(\lambda)A_j(\lambda).
$$
Substituting (\ref{02.7}) and (\ref{02.8}) into the last formula, we arrive at
\begin{equation}\label{02.10}
D_\nu^0(\lambda)=\sigma_m S_\nu'(1,\lambda)\prod_{l=2}^m v_{0,l}(\lambda) +\sigma_m S_\nu(1,\lambda)\sum_{j=2}^m v_{1,j}(\lambda)\prod_{{l\ne
j}\atop{l=2}}^m v_{0,l}(\lambda), \quad \nu=0,1.
\end{equation}
Comparing (\ref{02.10}) and (\ref{02.9}) with (\ref{02.4}) and (\ref{02.5}), respectively, we observe that the functions $D_\nu^0(\lambda)$
and $\Delta_\nu^0(\lambda)$ as well as $D_\nu(\lambda)$ and $\Delta_\nu(\lambda)$ differ only by one and the same sign $\sigma_m,$ which
finishes the proof.  $\hfill\Box$

\medskip
In particular, this lemma implies that the function $\Delta_\nu^0(\lambda)$ determined by formula (\ref{02.4}) is the characteristic
determinant of the problem ${\cal B}_\nu(0,H),$ while the function
$$
\Delta_\nu^{00}(\lambda)=S_\nu'(1,\lambda) (\cos\rho)^{m-1} +(1-m)S_\nu(1,\lambda)(\cos\rho)^{m-2}\rho\sin\rho
$$
is the characteristic determinant of ${\cal B}_\nu(0,0).$ Moreover, we have the asymptotics
$$
\Delta_\nu(\lambda)=\Delta_\nu^{00}(\lambda) +O\left(\rho^{\nu-1}e^{m|{\rm Im}\,\rho|}\right), \quad \lambda\to\infty, \quad \nu=0,1,
$$
while $\rho^{-\nu}\Delta_\nu^{00}(\rho^2)$ are sine-type functions of the exponential type $m$ in the $\rho$-plane. Hence, using Rouch\'e's
theorem as, e.g., in the proof of Theorem~4 in \cite{But22}, one can immediately conclude that zeros of $\Delta_\nu^{00}(\lambda)$ form a
principal part of the numbers $\lambda_{n,\nu}$ as $n\to\infty.$

Thus, analogously to Theorem~1.1.4 in \cite{novabook}, one can prove the following assertion.

\medskip
{\bf Lemma 2. }{\it The functions $\Delta_0(\lambda)$ and $\Delta_1(\lambda)$ are determined by their zeros uniquely. Specifically, the
representations
\begin{equation}\label{02.12}
\Delta_0(\lambda)=\prod_{n=1}^\infty\frac{\lambda_{n,0}-\lambda}{\lambda_{n,0}^0}, \quad
\Delta_1(\lambda)=m(\lambda_{1,1}-\lambda)\prod_{n=2}^\infty\frac{\lambda_{n,1}-\lambda}{\lambda_{n,1}^0}
\end{equation}
take place, where $\{\lambda_{n,\nu}^0\}_{n\in{\mathbb N}}$ are zeros of the function $\Delta_\nu^{00}(\lambda).$}
\\

{\large\bf 3. Solution of the inverse problem}
\\

Before proceeding directly to the proof of Theorem~1, we fulfil some preparatory work. Using Rouch\'e's theorem, one can prove the following
assertion.

\medskip
{\bf Lemma 3. }{\it For $j=\overline{2,m},$ zeros $\{\xi_{n,j}\}_{n\in{\mathbb N}}$ of the function $v_{0,j}(\lambda)$ have the asymptotics
\begin{equation}\label{3.1}
\xi_{n,j}=\eta_{n,j}^2, \quad \eta_{n,j}=\pi\Big(n-\frac12\Big) +\frac{H_j}{\pi n} +O\Big(\frac1{n^2}\Big), \quad n\to\infty.
\end{equation}
}

Without loss of generally, we assume that all multiple zeros $\xi_{n,j}$ are neighboring, i.e.
$$
\xi_{n,j}=\xi_{n+1,j}=\ldots=\xi_{n+m_{n,j}-1,j},
$$
where $m_{n,j}$ is the multiplicity of $\xi_{n,j}.$ Thus, the set
$$
{\cal S}_j:=\{n:\xi_{n,j}\ne \xi_{n-1,j},\, n-1\in{\mathbb N}\}\cup\{1\}
$$
indexes the zeros $\xi_{n,j}$ ignoring their multiplicities.

Similarly to Proposition~1 in \cite{BK19}, one can prove the following assertion, which remains valid for arbitrary sequences
$\{\xi_{n,j}\}_{n\in{\mathbb N}}$ of the form $\xi_{n,j}=\pi^2(n-1/2+\varkappa_n)^2$ with $\{\varkappa_n\}\in l_2.$

\medskip
{\bf Lemma 4. }{\it For any fixed $\nu\in\{0,1\}$ and $j=\overline{2,m},$ the functional system $\{s_{n,j,\nu}(x)\}_{n\in{\mathbb N}}$ forms
a Riesz basis in $L_2(0,1),$ where
\begin{equation}\label{3.2}
s_{k+s,j,\nu}(x)=k^{1-\nu} \frac{d^s}{d\lambda^s}S_\nu(x,\lambda)\Big|_{\lambda=\xi_{k,j}}, \quad k\in{\cal S}_j, \quad
s=\overline{0,m_{k,j}-1},
\end{equation}
while $S_\nu(x,\lambda)$ are determined in (\ref{02.0}).}

\medskip
For each $j=\overline{2,m}$ and $k\in{\cal S}_j,$ differentiating $p=\overline{0,m_{k,j}-1}$ times representation (\ref{02.5}) and
substituting $\lambda=\xi_{k,j}$ into the obtained relations, we arrive at the linear systems
\begin{equation}\label{3.3}
\gamma_{k+p,j,\nu} =\sum_{s=0}^p \Big({p\atop s}\Big) \alpha_{k+p-s,j,\nu} \beta_{k+s,j,\nu}, \quad p=\overline{0,m_{k,j}-1},\qquad k\in{\cal
S}_j, \quad \nu=0,1,
\end{equation}
where
\begin{equation}\label{3.3-1}
\gamma_{k+s,j,\nu}:=k^{1-\nu}(\Delta_\nu-\Delta_\nu^0)^{(s)}(\xi_{k,j}), \quad \alpha_{k+s,j,\nu} :=\frac{d^s}{d\lambda^s}\prod_{{l\ne
j}\atop{l=2}}^m v_{0,l}(\lambda) \Big|_{\lambda=\xi_{k,j}},
\end{equation}
\begin{equation}\label{3.4}
\beta_{k+s,j,\nu}:=k^{1-\nu}\frac{d^s}{d\lambda^s}V_j(Q_{\nu,j}(\,\cdot\,,\lambda))\Big|_{\lambda=\xi_{k,j}}, \quad k\in{\cal S}_j, \quad
s=\overline{0,m_{k,j}-1}.
\end{equation}

\medskip
{\bf Lemma 5. }{\it For $j=\overline{2,m},$ the following representations hold:
\begin{equation}\label{3.5}
V_j(Q_{0,j}(\,\cdot\,,\lambda)) =\frac{\omega_j}{2\rho} \sin\rho(2-a) +\int_0^{2-a} u_{0,j}(x) \frac{\sin\rho x}\rho\,dx,
\end{equation}
\begin{equation}\label{3.6}
V_j(Q_{1,j}(\,\cdot\,,\lambda)) = \frac{\omega_j}2\cos\rho(2-a) +H_j\omega_j\frac{\sin\rho(2-a)}{\rho} +\int_0^{2-a} u_{1,j}(x) \cos\rho
x\,dx,
\end{equation}
where
\begin{equation}\label{3.6-1}
\omega_j=\int_{a-1}^1 q_j(x)\,dx, \quad j=\overline{2,m},
\end{equation}
and
\begin{equation}\label{3.7}
u_{\nu,j}(x)=\frac14\Big(p_j\Big(\frac{a+x}2\Big) -(-1)^\nu p_j\Big(\frac{a-x}2\Big)\Big), \quad 0<x<2-a, \quad \nu=0,1,
\end{equation}
while
\begin{equation}\label{3.8}
p_j(x)=q_j(x)-2H_j\int_x^1 q_j(t)\,dt, \quad a-1<x<1.
\end{equation}
}

{\it Proof.} In accordance with (\ref{02.0}) and (\ref{02.1}), we have
$$
Q_{0,j}(x,\lambda)=\int_{a-1}^x \frac{\sin\rho(x-t)}\rho q_j(t) \frac{\sin\rho(t-a+1)}{\rho}\,dt \qquad\qquad\qquad\qquad\qquad\qquad
$$
\begin{equation}\label{3.9}
\qquad =\frac1{2\rho^2}\int_{a-1}^x\ q_j(t) \Big(\cos\rho(x-2t+a-1) -\cos\rho(x-a+1)\Big)\,dt,
\end{equation}
$$
Q_{1,j}(x,\lambda)=\int_{a-1}^x \frac{\sin\rho(x-t)}\rho q_j(t) \cos\rho(t-a+1)\,dt \qquad\qquad\qquad\qquad\qquad\qquad
$$
\begin{equation}\label{3.10}
\qquad =\frac1{2\rho}\int_{a-1}^x\ q_j(t) \Big(\sin\rho(x-2t+a-1) +\sin\rho(x-a+1)\Big)\,dt.
\end{equation}
By the differentiation, we obtain
$$
Q_{0,j}'(x,\lambda) =\frac1{2\rho}\int_{a-1}^x\ q_j(t) \Big(\sin\rho(x-a+1) -\sin\rho(x-2t+a-1)\Big)\,dt,
$$
$$
Q_{1,j}'(x,\lambda) =\frac12\int_{a-1}^x\ q_j(t) \Big(\cos\rho(x-a+1) +\cos\rho(x-2t+a-1)\Big)\,dt.
$$
Integrating by parts in (\ref{3.9}) and (\ref{3.10}), we arrive at
$$
Q_{0,j}(x,\lambda)=\frac1{\rho}\int_{a-1}^x \sin\rho(x-2t+a-1)\,dt\int_t^x q_j(\tau)\,d\tau,
$$
$$
Q_{1,j}(x,\lambda)=\frac{\sin\rho(x-a+1)}{\rho} \int_{a-1}^x q_j(t)\,dt -\int_{a-1}^x \cos\rho(x-2t+a-1)\,dt\int_t^x q_j(\tau)\,d\tau.
$$
Substituting $x=1$ into the last four formulae and using the designation (\ref{3.6-1}), we get
$$
Q_{0,j}'(1,\lambda) =\frac{\omega_j}{2\rho}\sin\rho(2-a) -\frac1{2\rho}\int_{a-1}^1 q_j(x)\sin\rho(a-2x)\,dx,
$$
$$
Q_{1,j}'(1,\lambda) =\frac{\omega_j}2\cos\rho(2-a) +\frac12\int_{a-1}^1 q_j(x)\cos\rho(a-2x)\,dx,
$$
$$
Q_{0,j}(1,\lambda) =\frac1{\rho}\int_{a-1}^1\sin\rho(a-2x)\,dx\int_x^1 q_j(t)\,dt,
$$
$$
Q_{1,j}(1,\lambda) = \frac{\omega_j}{\rho}\sin\rho(2-a) -\int_{a-1}^1\cos\rho(a-2x)\,dx\int_x^1 q_j(t)\,dt,
$$
which along with (\ref{1.11}) for $\nu_j=1$ and (\ref{3.8}) leads to the representations
$$
V_j(Q_{0,j}(\,\cdot\,,\lambda)) =\frac{\omega_j}{2\rho} \sin\rho(2-a) -\frac1{2\rho}\int_{a-1}^1 p_j(x)\sin\rho(a-2x)\,dx,
$$
$$
V_j(Q_{1,j}(\,\cdot\,,\lambda)) = \frac{\omega_j}2\cos\rho(2-a) +H_j\omega_j\frac{\sin\rho(2-a)}{\rho} +\frac12\int_{a-1}^1
p_j(x)\cos\rho(a-2x)\,dx.
$$
Changing the integration variable and taking (\ref{3.7}) into account, we arrive at (\ref{3.5}) and (\ref{3.6}).  $\hfill\Box$

\medskip
By virtue of (\ref{3.5}) and (\ref{3.6}) along with (\ref{02.0}) and (\ref{3.2}), formula (\ref{3.4}) takes the form
\begin{equation}\label{3.11}
\beta_{n,j,\nu} =\frac{\omega_j}2 s_{n,j,\nu}(2-a) +\delta_{\nu,1}H_j\omega_j \tilde s_{n,j,0}(2-a) +\int_0^{2-a}
u_{\nu,j}(x)s_{n,j,\nu}(x)\,dx, \quad n\in{\mathbb N},
\end{equation}
where $\delta_{\nu,1}$ is the Kronecker delta and
$$
\tilde s_{k+s,j,0}(x)=\frac1k s_{k+s,j,0}(x), \quad k\in{\cal S}_0, \quad s=\overline{0,m_{k,0}-1}.
$$

\medskip
{\bf Proof of Theorem~1.} By virtue of Lemma~2, specification of the spectra along with the coefficients $H_j$ uniquely determines the
numbers $\alpha_{n,j,\nu}$ and $\gamma_{n,j,\nu}$ defined in (\ref{3.3-1}). Further, since all $H_j$ are distinct, any two functions
$v_{0,j}(\lambda)$ and $v_{0,l}(\lambda)$ defined in (\ref{02.1}) for $j,l=\overline{2,m}$ have no common zeros whenever $j\ne l,$ i.e.
$\{\xi_{n,j}\}_{n\in{\mathbb N}} \cap\{\xi_{n,l}\}_{n\in{\mathbb N}}=\emptyset.$ Hence, $\alpha_{k,j,\nu}\ne0$ for all $k\in{\cal S}_j,$
$j=\overline{2,m}$ and $\nu=0,1.$ Therefore, each triangular system of linear algebraic equations (\ref{3.3}) always has the unique solution
$\beta_{k+s,j,\nu},$ $s=\overline{0,m_{k,j}-1}.$ Thus, all sequences $\{\beta_{n,j,\nu}\}_{n\in{\mathbb N}}$ are uniquely determined by
specifying the spectra of the problems ${\cal B}_0$ and ${\cal B}_1.$

Further, in accordance with (\ref{02.0}), (\ref{3.1}) and (\ref{3.2}), formula (\ref{3.11}) implies
$$
\frac{2\eta_{n,j}\beta_{n,j,0}}n =\omega_j\sin\eta_{n,j}(2-a) +2\int_0^{2-a} u_{0,j}(x)\sin\eta_{n,j}x\,dx, \quad n\gg1.
$$
Hence, the numbers $\omega_j$ are uniquely determined too by the formulae
\begin{equation}\label{3.12}
\omega_j=\lim_{k\to\infty}\frac{2\eta_{n_k,j}\beta_{n,j,0}}{n_k\sin\eta_{n_k,j}(2-a)}, \quad j=\overline{2,m},
\end{equation}
where $\{n_k\}$ is an increasing sequence of natural numbers that are chosen so that
$$
|a_{n_k}|\ge c>0, \quad a_n:=\sin\pi\Big(n-\frac12\Big)(2-a),
$$
whose existence, in turn, follows from the non-convergence of the sequence $\{a_n\},$ which can be shown similarly to Lemma~3.3 in
\cite{SatShieh19}.

Thus, by virtue of Lemma~4, all functions $u_{\nu,j}(x)$ are determined by relations (\ref{3.11}) uniquely. Hence, according to (\ref{3.7}),
the functions $p_j(x)$ are uniquely determined by the formulae
\begin{equation}\label{3.13}
 p_j(x)=
\left\{\begin{array}{cl}
 2(u_{1,j}-u_{0,j})(a-2x), & \displaystyle x\in\Big(a-1,\frac{a}2\Big),\\[3mm]
 2(u_{1,j}+u_{0,j})(2x-a), & \displaystyle x\in\Big(\frac{a}2,1\Big),
\end{array}\right. \quad j=\overline{2,m}.
\end{equation}
Finally, all $q_j(x)$ are unique solutions of the Volterra integral equations (\ref{3.8}), i.e.
\begin{equation}\label{3.14}
q_j(x)=p_j(x)+ 2H_j\int_x^1 e^{2H_j(t-x)} p_j(t)\,dt, \quad a-1<x<1, \quad j=\overline{2,m},
\end{equation}
which finishes the proof.  $\hfill\Box$

\medskip
This proof is constructive and gives the following algorithm for solving the inverse problem.

\medskip
{\bf Algorithm 1.} Let the spectra $\{\lambda_{n,0}\}_{n\in{\mathbb N}}$ and $\{\lambda_{n,1}\}_{n\in{\mathbb N}}$ be given. Then:

(i) Calculate the characteristic functions $\Delta_0(\lambda)$ and $\Delta_1(\lambda)$ by the formulae in (\ref{02.12});

(ii) Construct the sequences $\{\alpha_{n,j,\nu}\}_{n\in{\mathbb N}}$ and $\{\gamma_{n,j,\nu}\}_{n\in{\mathbb N}}$ for $j=\overline{2,m}$ and
$\nu=0,1$ by the formulae in (\ref{3.3-1}) and using (\ref{02.0})--(\ref{02.4});

(iii) Find the sequences $\{\beta_{n,j,\nu}\}_{n\in{\mathbb N}}$ for $j=\overline{2,m}$ and $\nu=0,1$ via relations (\ref{3.3});

(iv) Calculate the numbers $\omega_j, \,j=\overline{2,m},$ using (\ref{3.12});

(v) In accordance with (\ref{3.11}) and Lemma~4, construct the functions $u_{\nu,j}(x)\in L_2(0,2-a),$ by the formula
$$
u_{\nu,j}(x)=\sum_{n=1}^\infty \Big(\beta_{n,j,\nu} -\frac{\omega_j}2 s_{n,j,\nu}(2-a) -\delta_{\nu,1}H_j\omega_j \tilde
s_{n,j,0}(2-a)\Big)s_{n,j,\nu}^*(x), \;\; \nu=0,1, \;\; j=\overline{2,m},
$$
where $\{s_{n,j,\nu}^*(x)\}_{n\in{\mathbb N}}$ is the basis that is biorthogonal to the basis $\{\overline{s_{n,j,\nu}(x)}\}_{n\in{\mathbb
N}};$

(vi) Finally, find the functions $p_j(x)$ by (\ref{3.13}) and calculate $q_j(x)$ by (\ref{3.14}).
\\

{\bf Funding.} This research was supported by Russian Science Foundation, Grant No. 22-21-00509, https://rscf.ru/project/22-21-00509/

\medskip
{\bf Acknowledgement.} The author is grateful to Maria Kuznetsova for reading the manu\-script and making valuable comments.

\end{document}